# GENERICITY, THE ARZHANTSEVA-OL'SHANSKII METHOD AND THE ISOMORPHISM PROBLEM FOR ONE-RELATOR GROUPS

ILYA KAPOVICH AND PAUL SCHUPP

ABSTRACT. We apply the method of Arzhantseva-Ol'shanskii to prove that for an exponentially generic (in the sense of Ol'shanskii) class of one-relator groups the isomorphism problem is solvable in at most exponential time.

This is obtained as a corollary of the more general result that for any fixed integers $m > 1, n > 0$ there is an exponentially generic class of $m$-generator $n$-relator groups where every group has only one Nielsen equivalence class of $m$-tuples generating non-free subgroups. We also prove that all groups in this class are co-Hopfian.

## 1. INTRODUCTION

The idea of genericity in the context of geometric group theory was introduced by Gromov [23] and formalized by Ol'shanskii [47] who rigorously proved that "most" finitely presented groups are word-hyperbolic (see also the work of Champetier [11] for the two-relator case). The importance of probabilistic considerations for the study of finitely presented groups is now becoming increasingly clear, for example in the work of Gromov on "random" groups [24] and uniform embeddability into Hilbert spaces. See also [32, 33] for results on generic and average-case complexity of group-theoretic decision problems. Since the book of Gromov [23] appeared, there have been quite a number of results about genericity of various group-theoretic properties [2, 3, 4, 6, 11, 12, 13, 14, 15, 60, 62]. The basic philosophy is that generic behavior is often the best and the simplest possible. The results of this paper are very much in this line.

We now recall the formal definition of genericity according to Ol'shanskii [6]. Note that this definition is different from the definition used by Ol'shanskii in [47]. Roughly speaking, for fixed $m$ and $n$, a property of $m$-generated $n$-related groups is "generic" if a "randomly" chosen such group satisfies this property.

**Convention 1.1.** Throughout this paper, given $m$, we fix a finite alphabet $A = \{a_1, \ldots, a_m\}$, which will be the set of generators of the group $G$ under consideration. Let $F$ denote the free group $F(A) = F(a_1, \ldots, a_m)$ on $A$. A word $w$ in the alphabet $A \cup A^{-1}$ will be called *reduced* if it is freely reduced, that is if $w$ does not contain subwords of the form $a_i a_i^{-1}$ or $a_i^{-1} a_i$. A reduced word $w$ is said to be *cyclically reduced* if all cyclic permutations of $w$ are reduced. The *length* of a word $w$, denoted $|w|$, is the number of letters in $w$.

**Definition 1.2.** [6] Let $m > 1$ and $n > 0$ be fixed integers and fix the alphabet $A = \{a_1, \ldots, a_m\}$.

Let $P$ be a property of group presentations on $A$ of the form

$$\langle a_1, \ldots, a_m | r_1, \ldots, r_n \rangle$$

where the $r_i$ are cyclically reduced nontrivial words in $F(A)$.

For any integer $t \geq 0$ let $N(m,n,t)$ be the number of *all* possible presentations of the form $\langle a_1, \ldots, a_m | r_1, \ldots, r_n \rangle$ where the $r_i$ are cyclically reduced nontrivial words from $F(A)$ and where $|r_i| \leq t$ for $i = 1, \ldots, n$. Let $N_P(m,n,t)$ be the number of presentations with these restrictions which define a group with property $P$.

We say that $P$ is $(m,n)$-*generic* if

$$\lim_{t \to \infty} \frac{N_P(m,n,t)}{N(m,n,t)} = 1.$$

If, moreover, there is $0 \leq c = c(m,n) < 1$ such that for all sufficiently large $t$ we have

$$1 - \frac{N_P(m,n,t)}{N(m,n,t)} \leq c^t,$$

we say that $P$ is *exponentially $(m,n)$-generic*.

For the case $n = 1$ (one-relator groups) exponential genericity in the sense of Ol'shanskii coincides with Gromov's notion of exponential genericity [23, 47].

Let $G$ be a group. Recall that the classical *elementary Nielsen moves* on a tuple $(g_1, \ldots, g_m) \in G^m$ are: replace some $g_i$ by $g_i^{-1}$, or interchange some $g_i$ and $g_j$, or replace some $g_i$ by $g_i g_j$ for some $i \neq j$. Two $m$-tuples $\tau, \tau' \in G^m$ are *Nielsen-equivalent in $G$* if there is a finite chain of elementary Nielsen moves taking $\tau$ to $\tau'$. Clearly Nielsen-equivalent tuples generate the same subgroup of $G$. Another important basic fact is that in a free group of finite rank $F = F(x_1, \ldots, x_m)$ (where $m > 0$) every generating $m$-tuple of $F$ is Nielsen-equivalent to $(x_1, \ldots, x_m)$. This implies that for an arbitrary group $G$ two $m$-tuples $(g_1, \ldots, g_m)$ and $(g'_1, \ldots, g'_m)$ are Nielsen-equivalent in $G$ if and only if there is $\alpha \in Aut(F(x_1, \ldots, x_m))$, $\alpha(x_i) = W_i(x_1, \ldots, x_m)$ such that for $i = 1, \ldots, n$ we have $g'_i = W_i(g_1, \ldots, g_m)$.

The following theorem states that for a generic class of one-relator groups the isomorphism problem is solvable in at most exponential time. Moreover, given a one-relator presentation we can decide in exponential time if the presentation belongs to this generic class.

**Theorem A.** *Let $m > 1$ be an integer. There exists an exponentially $(m, 1)$-generic class $P_m$ of one-relator presentations $\langle a_1, \ldots, a_m | r \rangle$ with the following properties:*

1. *Every group $G$ defined by a presentation from $P_m$ is torsion-free, one-ended and word-hyperbolic. Moreover, every subgroup of $G$ generated by at most $m - 1$ elements is free.*
2. *There is an algorithm which, given an arbitrary cyclically reduced word $r \in F(a_1, \ldots, a_m)$, decides in at most exponential time (in the length of $r$) whether or not a presentation $\langle a_1, \ldots, a_m | r \rangle$ belongs to $P_m$.*
3. *For any presentation $\langle a_1, \ldots, a_m | r \rangle$ from $P_m$, for the group $G = \langle a_1, \ldots, a_m | r \rangle$ any $m$-tuple, generating a non-free subgroup of $G$, is Nielsen-equivalent in $G$ to the $m$-tuple $(a_1, \ldots, a_m)$.*
4. *Suppose $G_1 = \langle a_1, \ldots, a_m | r_1 \rangle$ and $G_2 = \langle a_1, \ldots, a_m | r_2 \rangle$ are one-relator presentations, at least one of which is in $P_m$. Then $G_1$ is isomorphic to $G_2$ if and only if there is $\alpha \in Aut(F(a_1, \ldots, a_m))$ such that either $\alpha(r_1) = r_2$ or $\alpha(r_1) = r_2^{-1}$.*
5. *There is an algorithm taking at most exponential time (in the sum of the lengths of the relators) which, given two $m$-generator one-relator presentations with at least one in $P_m$, decides if they define isomorphic groups.*

Theorem A is obtained as a corollary of the following more general statement:

**Theorem B.** *Let $m > 1$ and $n > 0$ be integers. There exists an exponentially $(m, n)$-generic class $P_{m,n}$ of $m$-generator $n$-relator presentations*

$$\langle a_1, \ldots, a_m | r_1, \ldots, r_n \rangle$$

*with the following properties:*

1. *Every group defined by a presentation from $P_{m,n}$ is torsion-free, one-ended and word-hyperbolic. Moreover, every subgroup of $G$ generated by at most $m - 1$ elements is free.*
2. *There is an algorithm which, given an arbitrary $m$-tuple of cyclically reduced words $r_1, \ldots, r_n \in F(a_1, \ldots, a_m)$, decides in at most exponential time (in the sum of the length of $r_i$) whether or not a presentation $\langle a_1, \ldots, a_m | r_1, \ldots, r_n \rangle$ belongs to $P_{m,n}$.*
3. *For any presentation $\langle a_1, \ldots, a_m | r_1, \ldots, r_n \rangle$ from $P_{m,n}$, for the group $G = \langle a_1, \ldots, a_m | r_1, \ldots, r_n \rangle$ any $m$-tuple, generating a non-free subgroup of $G$, is Nielsen-equivalent in $G$ to the $m$-tuple $(a_1, \ldots, a_m)$.*

Theorem B implies that any automorphism of a generic group

$$G = \langle a_1, \ldots, a_m | r_1, \ldots, r_n \rangle$$

from the class $P_{m,n}$ is induced by an automorphism of $F(a_1, \ldots, a_m)$. Moreover, $G$ has only one $m$-generated non-free subgroup, namely $G$ itself. This means, for example, that if $H$ is another $m$-generated group and $\phi : H \to G$ is a homomorphism then either $\phi(H) = G$ or $\phi(H)$ is free. In particular if $G = H$ and $\phi : G \to G$ is an injective endomorphism, then $\phi$ must be onto (since $G$ is non-free), so that $G$ is co-Hopfian. It is an important result of Sela [56] that all one-ended torsion-free hyperbolic groups are co-Hopfian. All groups covered by Theorem B are of this sort, but here we establish co-Hopficity for a generic class of finitely presented groups without using the JSJ-decomposition theory or the $\mathbb{R}$-tree techniques of Sela.

Theorem A is a corollary of Theorem B:

*Proof of Theorem A.* Parts (1),(2),(3) of Theorem A follow directly from the corresponding parts of Theorem B.

We say that an $m$-generator one-relator group has the *Nielsen uniqueness property* if this group has only one Nielsen equivalence class of generating $m$-tuples.

**Claim.** For an $m$-generated one-relator group $G = \langle a_1, \ldots, a_m | r \rangle$ with the Nielsen equivalence property and for any other $m$-generated one-relator group $H = \langle a_1, \ldots, a_m | s \rangle$ the groups $G$ and $H$ are isomorphic if and only if there is an automorphism of $F(a_1, \ldots, a_m)$ taking $s$ to either $r$ or $r^{-1}$.

The last condition is clearly sufficient for $G$ and $H$ to be isomorphic.

Suppose now that $\phi : H \to G$ is an isomorphism. To avoid confusion, we re-label the generators of $H$ and call them $b_i$ so that

$$H = \langle b_1, \ldots, b_m | s(b_1, \ldots, b_m) \rangle.$$

Denote $c_i = \phi(b_i)$ for $i = 1, \ldots, m$. The group $G$ is identified with $H$ via $\phi$ and hence on the generators $c_i$ the group $G$ has the presentation

$$(\dagger) \qquad G = \langle c_1, \ldots, c_m | s(c_1, \ldots, c_m) = 1 \rangle.$$

Since $G$ has the Nielsen uniqueness property, the $m$-tuple $(c_1, \ldots, c_m)$ is Nielsen-equivalent to $(a_1, \ldots, a_m)$ in $G$. Hence these tuples are related by an automorphism of $F(a_1, \ldots, a_m)$. When we rewrite presentation $(\dagger)$ of $G$ in the generators $(a_1, \ldots, a_m)$ using this automorphism, we conclude that for some $\alpha \in Aut(F(a_1, \ldots, a_m))$ the elements $r(a_1, \ldots, a_m)$ and $\alpha(s(a_1, \ldots, a_m))$ have the same normal closure in $F(a_1, \ldots, a_m)$. By a classical result of Magnus [40] this implies that $\alpha(s)$ is conjugate to $r$ or $r^{-1}$ in $F(a_1, \ldots, a_m)$, yielding the Claim. This establishes part (4) of Theorem A.

To see that (5) holds recall that by a classical result of Whitehead [41] one can decide in exponential time (in the sum of the lengths of the two words) whether two elements of $F(a_1, \ldots, a_m)$ are related by an automorphism. Applying this algorithm to the pairs $(r, s)$ and $(r^{-1}, s)$ yields the required result. □

The above algorithm also allows one to test the isomorphism of an arbitrary $m$-generator one-relator group and a group from a family of $m$-generator one-relator groups with the Nielsen equivalence property.

In terms of the paper [32] by Kapovich, Myasnikov, Schupp and Shpilrain Theorem A implies that the isomorphism problem for the class of all $m$-generated one-relator group is solvable strongly generically in exponential time. We refer the reader to [32] for the precise explanation of this terminology. We stress, however, that the above statement does not say that the isomorphism problem for $m$-generated one-relator groups is solvable. We also wish to point out that parts (1) of Theorem A and Theorem B are due to Arzhantseva and Olshanskii [6, 2] and we only state them for completeness.

There are very few general positive results regarding solvability of the isomorphism problem in combinatorial group theory. A famous exception is the remarkable theorem of Sela [56] asserting that the isomorphism problem is solvable in the class of torsion-free word-hyperbolic groups. One-relator groups are generically hyperbolic and so Sela's theorem applies to generic one-relator groups. However, besides using much more elementary considerations, the present theorem provides an explicit complexity bound for the algorithm solving the isomorphism problem and provides an exponential algorithm to decide if a given one-relator presentation does or does not belong to the generic class from Theorem A. Over the class of all finite presentations, there is no algorithm which decides if the group presented is word-hyperbolic since hyperbolicity is a Markov property. There is, however, a partial algorithm which eventually recognizes if a finitely presented group is hyperbolic [48]. It is still possible that hyperbolicity is decidable for one-relator presentations and this is a very interesting open question (see [30]).

The isomorphism problem for the general class of one-relator groups is still open. Pride [51] proved that the isomorphism problem is solvable for two-generator one-relator groups with torsion and Rosenberger [53, 54] established the same result for "cyclically pinched one-relator groups", that is for one-relator groups of the form $F_1 *_C F_2$ where $F_1, F_2$ are free and $C$ is infinite cyclic (see also [16, 61] for related work). The proofs of both Rosenberger and Pride essentially established the uniqueness of the Nielsen-equivalence class of generating tuples. This question of uniqueness has a long history of its own. It was a conjecture of Magnus [39] that every $m$-generator one-relator one-ended group has exactly one Nielsen equivalence class of generating $m$-tuples. McCool and Pietrowski [43] provided the first counter-example, where the number of Nielsen classes is bigger than one but finite. Brunner [9] constructed a two-generator one-relator group with infinitely many Nielsen equivalence classes of generating pairs. Pride [49] gave an example of a small cancellation (and hence hyperbolic) two-generated one-relator group where not every automorphism is "freely induced" and hence there is more than one Nielsen equivalence class of generating pairs. Moreover, it is possible to construct a two-generated infinitely related small cancellation group with infinitely many Nielsen equivalence classes of generating pairs [36]. Theorem A shows that, despite these important counterexamples, the conjecture of Magnus holds for an exponentially generic class of one-relator groups.

Other previously known related facts include a theorem of Gromov and Delzant [23, 19] that a torsion-free hyperbolic group has only finitely many (up-to conjugation) Nielsen-equivalence classes of pairs of elements generating one-ended subgroups. In particular, if $G$ is a two-generated torsion-free one-ended hyperbolic group then it has only finitely many Nielsen equivalence classes of generating pairs. We refer the reader to [21, 22, 28] and the further references therein for other results of a similar nature.

Moreover, Kapovich and Weidmann [38] proved that if $G$ is a torsion-free hyperbolic group where all $k$-generated subgroups are quasiconvex, then $G$ has only finitely many (up to conjugation) Nielsen-equivalence classes of $(k+1)$-tuples generating one-ended subgroups. Arzhantseva [3] proved, in particular, that in a generic $m$-generated $n$-related group all $(m-1)$-generated subgroups are quasiconvex (and free by the earlier result of Arzhantseva-Ol'shanskii [6]). Hence the theorem of Kapovich-Weidmann already implied that for a generic $m$-generated $n$-related group there are only finitely many Nielsen equivalence classes of generating $m$-tuples.

The main ingredient of our proofs is a powerful graph minimization technique introduced by Ol'shanskii and Arzhantseva [6]. We believe that this method is applicable in a variety of situations and deserves to be much more widely known and used. Their method is in many ways dual to the "perimeter reduction" technique of McCammond-Wise [8, 29, 42, 55]. Both methods allow one to study subgroups of non-free groups by means of Stallings subgroup graphs. Perimeter reduction reduces "what is missing" by filling in relator cycles. To be able to use this method, McCammond and Wise introduced the idea of a "distributive" small cancellation hypothesis which involves how occurrences of the generators are distributed among the relators. The Arzhantseva-Ol'shanskii minimization method reduces "what is present" by cutting out "large" parts of relator cycles. Arzhantseva and Ol'shanskii proposed a strong kind of "genericity" small cancellation condition adopted for groups with a given number of generators. We shall give a precise definition of this condition later.

The minimization technique was introduced by Arzhantseva and Ol'shanskii [6], where they used it to prove that for a generic $m$-generated $n$-related group all $(m-1)$-generated subgroups are free. It was later deployed by Arzhantseva [2, 3, 4] to study generic properties of finitely presented groups and subgroup properties of word-hyperbolic groups. (Results about groups where all subgroups with a bounded number of generators are free, were also obtained by Bumagina [10]). An illustration of the two different methods is provided by the work of Schupp [55] (using perimeter reduction) on subgroup separability of Coxeter groups and the paper of Kapovich and Schupp [34] using minimization with a "standard" small cancellation hypothesis to obtain results regarding freeness and quasiconvexity of subgroups of Coxeter groups, Artin groups and one-relator groups with torsion.

We now briefly describe the original idea of Arzhantseva-Olshanskii from [6]. Let $G$ be a small cancellation group with "good" genericity properties and with a finite generating set $A = \{a_1, \ldots, a_m\}$. If $\Gamma$ is a finite graph with edges labeled by letters from $A^{\pm 1}$ and base-vertex $v_0$, there is a canonical map $\phi : \pi_1(\Gamma, v_0) \to G$ sending a loop at $v_0$ to the element of $G$ represented by its label. Every $(m-1)$-generated subgroup $H$ of $G$ can be represented as $H = image(\phi)$ for at least one such graph $\Gamma$ with $rank(\pi_1(\Gamma)) \leq m - 1$ (for example, a wedge of circles labeled by the generators of $H$). Among all such $\Gamma$ representing $H$ with $rank(\pi_1(\Gamma)) \leq m - 1$ choose $\Gamma$ with with minimal *complexity* in some appropriate sense. The small cancellation assumption on $G$ then implies that either $\phi : \pi_1(\Gamma, v_0) \to G$ is injective and hence $H$ is free (the desired conclusion) or there is a reduced path in $\Gamma$ labeled by a large portion $w$ of a defining relator. Analysis of how this $w$ is subdivided by maximal arcs of $\Gamma$ allows Arzhantseva and Ol'shanskii to then perform a surgery trick on $\Gamma$ preserving the rank of its fundamental group as well as the subgroup $H = \phi(\pi_1(\Gamma, v_0)) \leq G$ but reducing the complexity of graph. This yields a contradiction with the minimal choice of $\Gamma$. We refer the reader to the proof of Lemma 5.1 below for a more detailed discussion on the Arzhantseva-Ol'shanskii method. In the present paper we make use of the fact [5] that both homotopy equivalence folding moves and the Arzhantseva-Ol'shanskii surgery trick translate into chains of Nielsen moves on the level of the generating sets $\phi(S)$ of $H$ where $S$ is a free basis of $\pi_1(\Gamma, v_0)$.

The authors are very grateful to Alexander Ol'shanskii and Goulnara Arzhantseva for useful comments on the preliminary versions of the paper.

2. Representing subgroups by labeled graphs

**Convention 2.1.** From now on, unless specified otherwise, we fix $m > 1$, $F = F(A) = F(a_1, \ldots, a_m)$ and a group $G = F(A)/N$, where $N$ is a normal subgroup of $F$.

Following the approach of Stallings [58], we use labeled graphs to study finitely generated subgroups of quotients of a free group.

**Definition 2.2.** An *F-graph* $\Gamma$ consists of an underlying oriented graph where every edge $e$ is labeled by a letter $\mu(e) \in A \cup A^{-1}$ in such a way that $\mu(e^{-1}) = \mu(e)^{-1}$ for every edge $e$ of $\Gamma$. We allow multiple edges between vertices as well as edges which are loops.

An $F$-graph $\Gamma$ is said to be *non-folded* if there exists a vertex $v$ and two distinct edges $e_1, e_2$ with origin $v$ such that $\mu(e_1) = \mu(e_2)$. Otherwise $\Gamma$ is said to be *folded*.

Every edge-path $p$ in $\Gamma$ has a label which is a word in $A \cup A^{-1}$. We shall denote this label by $\mu(p)$. The number of edges in $p$ will be called the *length* of $p$ and denoted $|p|$. A path $p$ in an $F$-graph $\Gamma$ is said to be *reduced* if it does not contains subpaths of the form $e, e^{-1}$ where $e$ is an edge of $\Gamma$. For a path $p$ (e.g. an edge) we will denote the initial vertex of $p$ by $o(p)$ and the terminal vertex of $p$ by $t(p)$.

The following statement is obvious:

**Lemma 2.3.** *Let $\Gamma$ be an $F$-graph. Then $\Gamma$ is folded if and only if the label of any reduced path in $\Gamma$ is a reduced word.*

Recall the definition of the classical folding move:

**Definition 2.4** (Fold)**.** Let $\Gamma$ be an $F$-graph. Suppose $e_1 \neq e_2$ are distinct edges of $\Gamma$ with common initial vertex $x = o(e_1) = o(e_2)$ and with the same label $a = \mu(e_1) = \mu(e_2) \in A \cup A^{-1}$. We fold the two edges $e_1, e_2$ into a single edge $e$ labeled $a$. The resulting $F$-graph $\Gamma'$ is said to be obtained from $\Gamma$ by a *fold*.

The following statement immediately follows from the definitions, exactly as in [58]:

**Proposition 2.5.** *Let $\Gamma$ be a connected graph and suppose that $\Gamma'$ is obtained from $\Gamma$ by a fold. Then the Euler characteristic of $\Gamma'$ is no less than that of $\Gamma$, that is $rank(\pi_1(\Gamma')) \leq rank(\pi_1(\Gamma))$.*

In addition to foldings, we need the following two transformations of labeled graphs introduced by Olshanskii and Arzhantseva [6].

Recall that we are working with a fixed presentation of a quotient $G$ of $F$.

**Definition 2.6** (Completing a relator cycle: move $M1$)**.** Let $p$ be path in an $F$-graph $\Gamma$ with initial vertex $x$, terminal vertex $y$ and label $\mu(p) = v$. Suppose that $v'$ is a reduced word such that $v = v'$ in $G$. We modify $\Gamma$ by attaching a new edge-path going from $x$ to $y$ labeled by the word $v'$.

We explicitly state the definition of the inverse move of $M1$:

**Definition 2.7** (Removing an arc from a relator cycle: move $M2$)**.** Let $p$ be a simple edge-path in a labeled $F$-graph $\Gamma$ with an initial vertex $x$, a terminal vertex $y$ and label $\mu(p) = v$ and suppose that $p$ is contained in a *maximal arc of* $\Gamma$ (an *arc* is a simple path, possibly closed, where every intermediate vertex of the path has degree two in $\Gamma$). Suppose there exists a path $p'$ in $\Gamma$ from $x$ to $y$ with label $\mu(p') = v'$ such that $p'$ and $p^{\pm 1}$ have no common edges and such that $v =_G v'$ in $G$.

We modify $\Gamma$ by removing all the edges of $p$ from $\Gamma$ while keeping the vertices $x$ and $y$.

Note that $M1$ decreases the Euler characteristic by one and that $M2$ increases the Euler characteristic by one.

**Definition 2.8** (Removing a degree-one vertex: move $R$)**.** Suppose $e$ is an edge of an $F$-graph $\Gamma$ such that the vertex $t(e)$ has degree one in $\Gamma$. We remove the edge $e$ and the vertex $t(e)$ from $\Gamma$.

**Definition 2.9** (Combination Arzhantseva-Ol'shanskii move: move $AO$). Suppose $\Gamma$ is a connected $F$-graph. Let $p_1 p' p_2$ be a reduced path in $\Gamma$ such that $p'$ is a path contained in a maximal arc of $\Gamma$ and the paths $p_1$, $p_2$ do not overlap $p'$. Let $u_1$ and $u_2$ be respectively the labels of $p_1$ and $p_2$ and let $u$ be the label of $p'$. Suppose $y$ is a reduced word such that $u_1 u u_2 y = 1$ in $G$. Perform a combination of $M1$ and $M2$ by first attaching to $\Gamma$ a new edge-path $f$ labeled $y$ from $t(p_2)$ to $o(p_1)$ and then removing the arc $p'$. The resulting graph $\Gamma'$ is said to be obtained from $\Gamma$ by a *move of type AO*.

**Proposition-Definition 2.10.** Let $\Gamma$ be a connected $F$-graph with a base-vertex $x_0$. Then the labeling of paths gives rise to a homomorphism
$$\phi : \pi_1(\Gamma, x_0) \to G$$
such that for every path $p$ from $x_0$ to $x_0$ we have $\phi([p]) =_G \mu(p)$ in $G$ (where $[p]$ stands for the equivalence class of $p$ in $\pi_1(\Gamma, x_0)$). In this case we will say that $H = \phi(\pi_1(\Gamma, x_0)) \leq G$ is the *subgroup represented by* $(\Gamma, x_0)$.

Moreover, the following holds:
(1) If $\Gamma$ is finite then $image(\phi)$ is finitely generated. If, in addition, $\Gamma$ has Euler characteristic $1 - k$, then the free group $\pi_1(\Gamma, x_0)$ has rank $k$ and hence $\phi(\pi_1(\Gamma, x_0))$ can be generated by $k$ elements.
(2) If $x_1$ is another vertex of $\Gamma$ then the pairs $(\Gamma, x_0)$ and $(\Gamma, x_1)$ define conjugate subgroups of $G$.
(3) Every finitely generated subgroup of $G$ can be represented in this fashion for some finite connected $\Gamma$. Moreover, if $H \leq G$ is $k$-generated, then $H$ can be represented by a connected $F$-graph of Euler characteristic $\geq 1 - k$.

The following simple fact plays an important role in our approach.

**Proposition 2.11.** *Let $\Gamma$ be an $F$-graph with a base-vertex $x_0$. Suppose $\Gamma'$ is obtained from $\Gamma$ by a finite sequence of folds and as well as moves $M1, M2$, so that $x_0'$ is the image of $x_0$ in $\Gamma'$. Then the pairs $(\Gamma, x_0)$ and $(\Gamma', x_0')$ define the same subgroup of $G$.*

*Moreover, if $\Gamma'$ is obtained from $\Gamma$ by removing a degree-one vertex move $R$, then for any vertices $x$ of $\Gamma$ and $x'$ of $\Gamma'$ the pairs $(\Gamma, x)$ and $(\Gamma', x')$ define conjugate subgroups of $G$.*

## 3. Nielsen equivalence

We wish to examine more closely the effects of the moves $AO, M1, M2$ on Nielsen-equivalence classes of generating tuples for the subgroups defined by labeled graphs.

The following simple but crucial lemma is due to Arzhantseva (Lemma 2 of [5]). The proof is a straightforward corollary of the definitions and relies on the fact that any two free bases of a free group of finite rank are Nielsen-equivalent. This applies in particular to the fundamental group of a finite graph.

**Lemma 3.1.** [5] *Let $\Gamma$ be a connected $F$-graph with a base-vertex $x_0$. Let $\Gamma'$ be obtained from $\Gamma$ by a fold preserving the Euler characteristic or by a move AO or by removing a degree-one vertex different from $x_0$. Denote the image of $x_0$ by $x_0'$. Let $\phi : \pi_1(\Gamma, x_0) \to G$ and $\phi' : \pi_1(\Gamma', x_0') \to G$ be the label maps for $\Gamma$ and $\Gamma'$ respectively. Thus $\mathrm{image}(\phi) = \mathrm{image}(\phi) = H \leq G$.*

*Then for any $m$-tuples $\tau$ and $\tau'$ freely generating $\pi_1(\Gamma, x_0)$ and $\pi_1(\Gamma', x_0')$ respectively, the $m$-tuples $\phi(\tau)$ and $\phi(\tau')$ are Nielsen-equivalent in $G$.*

*Moreover, if $x_0$ has degree one in $\Gamma$ and $\Gamma'$ is obtained from $\Gamma$ by an $R$-move removing $x_0$ then for any vertex $x_0'$ of $\Gamma'$ and any $m$-tuples $\tau$ and $\tau'$ as above, the $m$-tuple $\phi(\tau)$ is conjugate to an $m$-tuple that is Nielsen-equivalent to $\phi(\tau')$ in $G$.*

## 4. The Genericity condition

We fix $m > 1$ and thus $A = \{a_1, \ldots, a_m\}$ and $F = F(A)$. We now give the definitions needed for the genericity condition.

**Definition 4.1.** [6] Let $0 < \mu \leq 1$ be a real number. A reduced word $w$ in $F(A) = F(a_1, \ldots, a_m)$ of length $l > 0$ is called *$\mu$-readable* if there exists a connected $F$-graph $\Gamma$ where every edge is labeled by a letter of $A$ such that:

  (1) The number of edges in $\Gamma$ is at most $\mu l$.
  (2) The free group $\pi_1(\Gamma)$ has rank at most $m - 1$.
  (3) There is a reduced path in $\Gamma$ with label $w$.

**Definition 4.2.** [2] Let $0 < \mu \leq 1$ be a real number and let $L > 1$ be an integer. A reduced word $w$ in $F(A)$ of length $l > 0$ is called *$(\mu, L)$-readable* if there exists a connected $F$-graph $\Gamma$ such that:

  (1) The number of edges in $\Gamma$ is at most $\mu l$.
  (2) The free group $\pi_1(\Gamma)$ has rank at most $L$.
  (3) There is a path in $\Gamma$ with label $w$.
  (4) The graph $\Gamma$ has at least one vertex of degree $< 2m$.

Clearly, for fixed integer $L$ and a rational $\mu$ the problem of deciding if an arbitrary $w \in F(A)$ is $\mu$-readable or $(\mu, L)$-readable, is decidable in time exponential in the length of $w$. Although the definition may seem technical, how the concept of $(\mu, L)$-readability leads to constructing exponentially generic families of presentations is actually simple. The only folded subgroup graph of rank $m$ for a set of generators of the free group $F$ of rank $m$ is the bouquet of $m$ loops at a single vertex where the loops are labeled by the generators $a_i$ of $F$. All words of a given length $l$ can be obtained as the label of a reduced path in this graph. Compare this case with the other graph of rank $m = 3$ which has a vertex of degree less than $2m$, as illustrated in Figure 1 for $F = F(a, b, c)$. It is easy to see that the number of words which are readable as labels of reduced paths of length $l$ in the second graph divided by the number of all possible reduced words of length $l$ is going to $0$ exponentially fast.

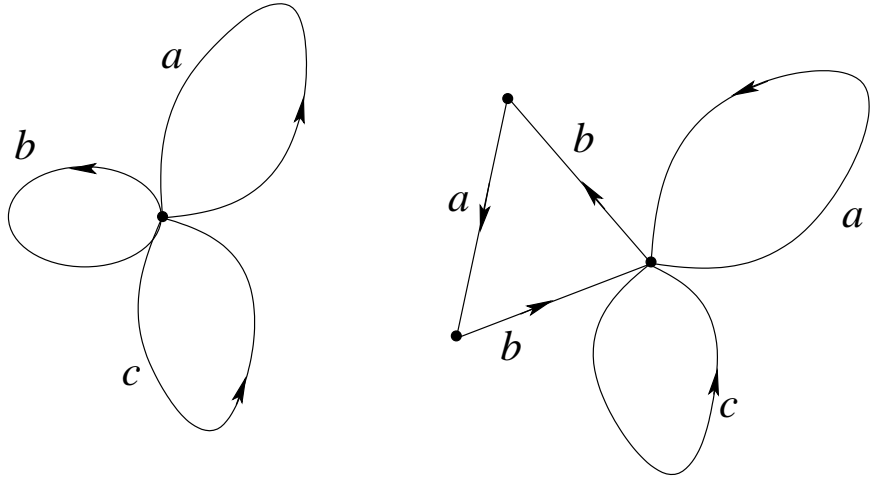

FIGURE 1. Illustrating the "readability" condition

We can now define an appropriate small cancellation condition to take advantage of such a situation.

**Definition 4.3.** Let $0 < \mu \leq 1$ be a real number, let $n > 0$ and $L \geq m$ be integers and let $0 < \lambda$ be a real number such that
$$\lambda \leq \frac{\mu}{15L + 3\mu} \leq \frac{\mu}{15m + 3\mu} < 1/6.$$
We will say that a tuple of nontrivial cyclically reduced words $(r_1, \ldots, r_n)$ in $F(A)$ satisfies the $(\lambda, \mu, L)$-condition if:

(1) The tuple $(r_1, \ldots, r_n)$ satisfies the $C'(\lambda)$ small cancellation condition.
(2) The words $r_i$ are not proper powers in $F(A)$.
(3) If $w$ is a subword of a cyclic permutation of some $r_i$ and $|w| \geq |r_i|/2$ then $w$ is not $(\mu, L)$-readable and not $\mu$-readable.

We refer the reader to [41] for the basic definitions and background information regarding small cancellation theory.

Again, it is easy to see that if we fix an integer $L$ and rational numbers $\lambda, \mu$ as above then for an arbitrary $n$-tuple $(r_1, \ldots, r_n)$ we can decide whether or not the tuple satisfies the $(\lambda, \mu, L)$-condition in time at most exponential in $\sum |r_i|$. The following condition is the intersection of the $(\lambda, \mu)$-condition from [6] and what was called $(\lambda, \mu, L)$-condition in [2].

**Definition 4.4.** Let $m > 1$, $n > 0$ and $F = F(a_1, \ldots, a_m)$. Let $\lambda, \mu, L$ be as in Definition 4.3. We define the class $P(\lambda, \mu, L, m, n)$ of finite group presentations on $m$ generators with $n$ defining relators as the collection of all presentations of the form
$$\langle a_1, \ldots, a_m | r_1, \ldots r_n \rangle$$
where the $n$-tuple $(r_1, \ldots, r_n)$ satisfies the $(\lambda, \mu, L)$-condition.

Arzhantseva and Ol'shanskii [6, 2] prove the following theorem.

**Theorem 4.5.** *The condition $P(\lambda, \mu, L, m, n)$ is exponentially $(m,n)$-generic.*

## 5. Proofs of the main results

The following lemma plays an important role in our argument. This statement is due to Arzhantseva [2] and is contained (although not stated explicitly) in section 4 of [2]. In fact, this lemma is a crucial step in the proof of Theorem 1 of [2].

**Lemma 5.1.** [2] *Let $G$ be a group given by a finite presentation*

$$G = \langle a_1, \ldots, a_m | r_1, \ldots, r_n \rangle$$

*belonging to class $P(\lambda, \mu, L, m, n)$.*

*Let $\Gamma$ be a connected folded F-graph with a base-vertex $x_0$ and with*

$$rank(\pi_1(\Gamma), x_0)) \le L.$$

*Suppose that $\Gamma$ has at least one vertex has degree $< 2m$ and no vertices of degree one, except possibly $x_0$.*

*Then either $\phi : \pi_1(\Gamma, x_0) \to G$ is injective (and hence $H = image(\phi)$ is free) or there exists an AO-move on $\Gamma$ that reduces the number of edges of $\Gamma$.*

*Proof.* For completeness, we will give a sketch of the argument, which first appeared, in essentially the same form, in the paper of Arzhantseva-Ol'shanskii [6]. If $\phi$ is not injective then there is a nontrivial reduced loop at $x_0$ in $\Gamma$ whose label is equal to 1 in $G$. Since $\Gamma$ is folded, the label of this loop is a reduced word. The small cancellation assumption $C'(\lambda)$ on the presentation of $G$ then implies that the label of this loop contains a subword $v$ which is also a subword of a cyclic permutation $r$ of some $r_i^{\pm 1}$ from the presentation of $G$ and such that $v$ is missing less than $3\lambda|r|$ letters of $r$. After possibly replacing the loop by its inverse, we may assume that $v = r_i$. Thus $v$ is the label of a reduced path $p$ in $\Gamma$ and $|v| > (1 - 3\lambda)|r|$. The assumption that $rank(\pi_1(\Gamma), x_0)) \le L$ implies that $\Gamma$ contains at most $3L - 1$ maximal arcs. These maximal arcs break up $p$ as a concatenation $p = p_1 p_2 \ldots p_k$ of $k$ subpaths with labels $v_1, \ldots, v_k$, where $p_2 \ldots p_{k-1}$ are maximal arcs. Thus $v = v_1 \ldots v_k$.

Suppose first that there is some $p_i$ with $|p_i| \ge 5\lambda|r|$. In this situation a case-by-case analysis of possible overlaps of $p_i$ with the other paths $p_j$ will show that some $AO$-type move reduces the number of edges.

We claim that there is a subpath of $p_i$ of length at least $3\lambda|r|$ that does not overlap the rest of the path $p$. Indeed, suppose first that $1 < i < k$, so that $p_i$ is a maximal arc of $\Gamma$. Since $|p_i| \ge 5\lambda|r|$, the $C'(\lambda)$-condition and the assumption that $r$ is not a proper power imply that $p_i \ne p_j^{\pm 1}$ for $j \ne i, 1 < j < k$ and hence $p_i$ does not overlap any of such $p_j$. Similarly, the overlap of $p_i$ with each of $p_1, p_k$ has length less than $\lambda|r|$. Therefore there is a subpath of $p_i$ of length at least $3\lambda|r|$ that does not overlap the rest of the path $p$, as claimed. Suppose now that $i = 1$ (the case $i = k$ is symmetric) and $|p_1| \geq 5\lambda|r|$. Then $p_1$ does not overlap any of $p_i$, $1 < i < k$ since otherwise we are in the previous case. The overlap between $p_1$ and $p_k$ has length less than $\lambda|r|$ because of the $C'(\lambda)$-condition and the fact that $r$ is not a proper power. This yields the claim as before.

Thus indeed there is a subpath $p'$ of $p_i$ such that $|p'| \geq 3\lambda|r|$ and that $p'$ does not overlap the rest of $p$. Then we can perform an $AO$-move on $\Gamma$ as follows: first add an arc from $t(p_k) = t(p)$ to $o(p_1) = o(p)$ labeled by the missing in $v$ part of the relator $r$ (of length $< 3\lambda|r|$) and then remove the arc $p'$. Clearly, this reduces the number of edges in $\Gamma$ and the conclusion of Lemma 5.1 holds.

Suppose now that for each $i = 1, \ldots, k$ we have $|p_i| < 5\lambda|r|$. Then the word $v$ is readable as a label of a path in some connected subgraph $\Gamma_1$ of $\Gamma$ and the rank of the fundamental group of $\Gamma_1$ is at most $L$: namely, take $\Gamma_1$ to be the union of all edges traversed by $p$, that is the edges contained in $\cup_{i=1}^{k} p_i$. The graph $\Gamma$ contains at least one vertex of degree $< 2n$ by assumption and hence so does $\Gamma_1$ (Lemma 2 of [2]). Moreover, $\Gamma$ has at most $3L - 1$ distinct maximal arcs and each $p_i$ is contained in a maximal arc of $\Gamma$. Therefore the number of edges in $\Gamma_1$ is $\leq 3L \cdot 5\lambda|r| \leq \mu(1 - 3\lambda)|r| \leq \mu|v|$ by the choice of $\lambda$ and $v$. Since $v$ is a subword of a defining relation $r$ with $|v| \geq |r|/2$, we get a contradiction with the assumption that the presentation of $G$ belongs to $P(\lambda, \mu, L, m, n)$.

$\square$

**Theorem 5.2.** *Let $L \geq m > 1, n > 0$ be integers and let $F = F(A) = F(a_1, \ldots, a_m)$. Let $0 < \mu \leq 1$ and $\lambda > 0$ be real numbers such that*

$$\lambda \leq \frac{\mu}{15L + 3\mu} \leq \frac{\mu}{15m + 3\mu} < 1/6.$$

*Let $G$ be a group defined by a presentation $G = \langle a_1, \ldots, a_m | r_1, \ldots r_n \rangle$ that belongs to the class $P(\lambda, \mu, L, m, n)$.*

*Then $G$ is torsion-free one-ended word-hyperbolic and every $(m - 1)$-generated subgroup of $G$ is free. Moreover any $m$-tuple, generating a non-free subgroup of $G$, is Nielsen-equivalent to the tuple $(a_1, \ldots, a_m)$ in $G$.*

*Proof.* The results of Arzhantseva [2] imply that $G$ is freely indecomposable, torsion-free, non-elementary, and word-hyperbolic. In particular, she proved that any group in class $P(\lambda, \mu, L, m, n)$ is not a free group.

We now come to the only place where we use the $\mu$-readability assumption. Namely, $G$ cannot be generated by fewer than $m$ elements since Arzhantseva-Ol'shanskii [6] proved that every $(m - 1)$-generated subgroup of $G$ is free. Hence $G$ must have rank exactly $m$.

Let $(w_1, \ldots, w_m)$ be an $m$-tuple of nontrivial freely reduced words in $A \cup A^{-1}$ generating a non-free subgroup $H$ of $G$. We will show that $(w_1, \ldots, w_m)$ is Nielsen-equivalent to $(a_1, \ldots, a_m)$ in $G$.

We may assume that $\sum |w_i| > m$ since otherwise $(w_1, \ldots, w_m)$ is a permutation of $(a_1^{\pm 1}, \ldots, a_m^{\pm 1})$ and there is nothing to prove. We will inductively define a sequence of finite connected $F$-graphs $\Gamma_i$ (for $i = 0, 1, \dots$) with the following properties:

a) For each $i$ the graph $\Gamma_i$ represents a conjugate of $H$, that is for the label-maps $\phi_i : \pi_1(\Gamma_i) \to G$ we have $image(\phi_i) = g_i H g_i^{-1}$ for some $g_i \in G$ (since $\Gamma_i$ represents a conjugate of $H$ and is connected, the choice of a base-vertex is irrelevant here).

b) For each $i$ the free group $\pi_1(\Gamma_i)$ has rank $m$.

c) For each $i$ the graph $\Gamma_{i+1}$ (if it is defined) has fewer edges than $\Gamma_i$.

d) Each $\Gamma_{i+1}$ is obtained from $\Gamma_i$ by a combination of foldings, moves of type $AO$, and removing vertices of degree one.

Note that $H$ has rank $m$ since otherwise $H$ can be generated by fewer than $m$ elements and hence is free by the properties of $G$.

Let $\Gamma_0$ be a wedge of $m$ loop-edges labeled by $w_1, \dots, w_m$. Note that $w_1, \dots, w_m$ generate a free group of rank $m$ in $F(A)$ since otherwise the group $H$ can be generated by fewer than $m$ elements.

We will now describe the inductive step for constructing $\Gamma_{i+1}$ from $\Gamma_i$.

Let $\Gamma'$ be obtained from $\Gamma_i$ by performing a sequence of foldings and removing degree-one vertices until no more such moves are possible. Thus $\Gamma'$ is folded, connected and has no degree-one vertices. The folding moves do not increase the rank of the fundamental group of an $F$-graph. Moreover, $rank(\pi_1(\Gamma')) = rank(\pi_1(\Gamma_i)) = m$ since $rank(\pi_1(\Gamma')) < m$ would imply that $H$ could be generated by fewer than $m$ elements.

If $\Gamma'$ is the wedge of $m$ loop-edges labeled $a_1, \dots, a_m$, we put $\Gamma_{i+1} = \Gamma'$ and terminate the process. Suppose not. Then the map $\phi : \pi_1(\Gamma') \to G$ (with $image(\phi)$ being a conjugate of $H$) is not injective since $H$ is not free. Note that there is at least one vertex of degree $< 2m$ in $\Gamma'$. Otherwise $\Gamma'$ defines a subgroup of finite index $k > 1$ in $F(a_1, \dots, a_m)$ and hence $rank(\pi_1(\Gamma')) = k(m-1) + 1 > m$ by the Schreier index formula, contrary to the fact that $rank(\pi_1(\Gamma')) = m$.

Lemma 5.1 then implies that there is a move of type $AO$ which, when applied to $\Gamma'$, produces a graph $\Gamma''$ with smaller total number of edges. Note that by definition an $AO$-move does not change the rank of the fundamental group of a graph. Put $\Gamma_{i+1} = \Gamma''$.

It is clear that conditions a)-d) above hold for $\Gamma_{i+1}$. This completes the description of the inductive step.

Since the total number of edges decreases in the above process, the sequence $\Gamma_0, \Gamma_1, \dots$ must terminate in a finite number of steps with the graph $\Gamma_q$ which is the wedge of $m$ loop-edges labeled $a_1, \dots, a_m$.

By applying Lemma 3.1 to the above sequence of graphs and looking at the free bases of $\pi_1(\Gamma_0)$ and $\pi_1(\Gamma_q)$ we conclude that $(w_1, \dots, w_m)$ is conjugate to a tuple that is Nielsen-equivalent to $(a_1, \dots, a_m)$ in $G$. Hence $(w_1, \dots, w_m)$ generates $G$. For $m$-tuples generating $G$ conjugation implies Nielsen-equivalence and hence $(w_1, \dots, w_m)$ is Nielsen-equivalent to the tuple $(a_1, \dots, a_m)$ in $G$, as required.

□

*Proof of Theorem B.* The results of Arzhantseva-Ol'shanskii (Lemma 4 of [6]) and Arzhantseva (Lemma 5 of [2]) imply that under the assumptions of Theorem 5.2 the class $P(\lambda, \mu, L, m, n)$ is exponentially $(m,n)$-generic (since it is the intersection of two exponentially generic classes). Hence Theorem 5.2 above, by choosing rational $\mu$ and $\lambda$, implies Theorem B from the introduction. □

Department of Mathematics, University of Illinois at Urbana-Champaign, 1409 West Green Street, Urbana, IL 61801, USA
*E-mail address*: `kapovich@math.uiuc.edu`

Department of Mathematics, University of Illinois at Urbana-Champaign, 1409 West Green Street, Urbana, IL 61801, USA
*E-mail address*: `schupp@math.uiuc.edu`